\documentclass{amsart}
\usepackage{amssymb,amsmath,amsthm}
\usepackage{mathrsfs}

\usepackage{lineno}
\usepackage{tikz}

\newcommand{\Sym}{ \mathfrak{S}}
\renewcommand{\ell}{l}
\newcommand{\x}{\mathbf{x}}
\newcommand{\Orb}{\mathscr{O}}
\newcommand{\pfp}[1]{N\left(#1!\right)}
\newcommand{\mfp}[1]{N\big(f(#1)!\big) }
\newcommand{\one}{\mathbf{1}}

\newcommand{\pp}[1]{N\left(#1\right)}

\newcommand{\MD}[1]{M\!D\left(#1\right)}

\newtheorem{theorem}{Theorem}

\newtheorem{observation}[theorem]{Observation}

\theoremstyle{definition}
\newtheorem{remark}[theorem]{Remark}
\newtheorem{definition}[theorem]{Definition}

\numberwithin{equation}{section}
\numberwithin{theorem}{section}

\begin{document}

\title{The combinatorics of\\ MacMahon's partial fractions}

\author{Andrew V. Sills}

\thanks{The author thanks the National Security Agency for partially supporting his research
program via grant H98230-14-1-0159 during 2014--2015, when this research project commenced.}

\email{ASills@georgiasouthern.edu}
\urladdr{http://home.dimacs.rutgers.edu/~asills}
\date{\today}

\keywords{partitions;  partition function; compositions; symmetric group}
\subjclass[2010]{Primary 05A17}

\maketitle

\centerline{Dedicated to my teacher, George E. Andrews, on the occasion of his 80th
birthday}

\begin{abstract}
MacMahon showed that the generating function for partitions into at most $k$ parts can be decomposed into a partial fractions-type sum indexed by the partitions of $k$.  In this present work, a generalization of MacMahon's result is given, which in turn provides a full combinatorial explanation. 
\end{abstract}

\section{Introduction}
\subsection{An excerpt from MacMahon's {\it Combinatory Analysis}}
In \emph{Combinatory Analysis}, vol. 2~\cite[p. 61ff]{M2}, P. A. MacMahon writes:
 \begin{quote}
   We commence by observing the two identities
   \begin{align*}
   \frac{1}{(1-x)(1-x^2)} &= \frac{1}{2} \frac{1}{(1-x)^2} + \frac 12 \frac{1}{(1-x^2)},\\
 \frac{1}{(1-x)(1-x^2)(1-x^3)} &= \frac{1}{6} \frac{1}{(1-x)^3} + \frac 12 \frac{1}{(1-x)(1-x^2)}+
 \frac 13 \frac{1}{1-x^3},
 \end{align*} 
 which we will also write in the illuminating notation so often employed:
 \begin{align*}
   \frac{1}{\mathbf{(1)(2)}}&= \frac{1}{2} \frac{1}{\mathbf{(1)^2}} + \frac 12 \frac{1}{\mathbf{(2)}},\\
 \frac{1}{\mathbf{(1)(2)(3)}} &= \frac{1}{6} \frac{1}{\mathbf{(1)^3}} + \frac 12 \frac{1}{\mathbf{(1)(2)}}+
 \frac 13 \frac{1}{\mathbf{(3)}}.
 \end{align*}
 [Elsewhere~\cite[p. 5]{M2} the ``very convenient notation" $\mathbf{(j)} = 1-x^j$ is defined and attributed to Cayley.]
 
   The observation leads to the conjecture that we are in the presence of partial fractions of a new and special kind.  We note that in the first identity we have a fraction corresponding to each of the partitions $( 1^2)$, $(2)$ of the number $2$ and in the second fractions corresponding to and derived from each of the partitions $(1^3)$, $(21)$, $(3)$ of the number $3$. \dots
   [In general] we find
   \begin{equation} \label{Meq}
     \frac{1}{ \mathbf{(1)(2)\cdots(i)}} =
     \sum \frac{1}{1^{p_1} . 2^{p_2} . 3^{p_3} \dots  p_1! p_2! p_3! \dots} 
     \frac{1}{\mathbf{ (1)}^{p_1} \mathbf{(2)}^{p_2} \mathbf{(3)}^{p_3}  \dots },
   \end{equation}
where $(1^{p_1} 2^{p_2} 3^{p_3} \dots)$ is a partition of $i$ and the summation is in regard to all partitions of $i$.   
  This remarkable result shows the decomposition of the generating function into as many fractions as the number $i$ possesses partitions.  The denominator of each fraction is directly derived from one of the partitions and is of degree $i$ in $x$.  The numerator does not involve $x$ and the coefficient is the easily calculable number
  \[ \frac{1}{1^{p_1} . 2^{p_2} . 3^{p_3} \dots  p_1! p_2! p_3! \dots}. \]
 \end{quote}
 
 \begin{remark}
 In~\cite[p. 209, Ex. 1]{A76} Andrews attributes~\eqref{Meq} to Cayley, and this attribution 
has been repeated by other authors in the literature.  However, the author 
has been unable to find~\eqref{Meq} anywhere in Cayley's works~\cite{C2}, and indeed 
MacMahon, in the chapter where he presents his ``partial fractions of a new and special 
kind"~\cite[Section VII, Chapter V]{M2} contrasts his results with those of Cayley several 
times.
 \end{remark}

\subsection{Some definitions and notation}
\subsubsection{Partitions and related objects}
It will be necessary to employ partitions and compositions of positive integers, sometimes allowing $0$'s as parts, and sometimes not.  Accordingly, we will formalize terminology via the following definitions.

\begin{definition}\label{defweakcomp}
A \emph{weak $k$-composition} $\gamma$ is a $k$-tuple of nonnegative integers 
$(\gamma_1, \gamma_2, \dots, \gamma_k)$.  Each $\gamma_i$ (even if $\gamma_i=0$) is called a \emph{part} of $\gamma$.  The \emph{weight} of $\gamma$, denoted $|\gamma|$, is $\sum_{i=1}^m \gamma_i$.
The \emph{length} of $\gamma$, denoted $\ell(\gamma)$, is the number of parts in $\gamma$. 
The \emph{frequency} (or \emph{multiplicity})  
of part $j$ in $\gamma$, denoted $f_j(\gamma)$ or simply $f_j$ when 
$\gamma$ is clear from context, is the number of times
that $j$ appears as a part in $\gamma$:
\[ f_j(\gamma) := \# \{ i : \gamma_i = j \}. \] 
The \emph{frequency sequence} associated with $\gamma$ is
\[ f(\gamma) := ( f_0(\gamma), f_1(\gamma), f_2(\gamma), f_3(\gamma), \dots ). \]
 The set of all weak $k$-compositions will be denoted $\mathscr{C}_k$.
\end{definition}

\begin{remark}
The author prefers to use the term \emph{frequency} and corresponding notation
$f_j$ over the term \emph{multiplicity} (with notation $m_j$) to be consistent with the
works of Andrews~\cite{A67,A72,A76}.
\end{remark}

\begin{definition}
A weak $k$-composition $w=(w_1,w_2,\dots,w_k)$ 
is a \emph{weak $k$-partition} if its parts occur in 
nonincreasing order
\[ w_1 \geq w_2 \geq \cdots \geq w_k. \]  The set of weak $k$-partitions of weight $n$ will be denoted $\mathscr{W}_{k}(n)$ and the cardinality of this set by $p_k(n)$.  
\end{definition}

\begin{definition}
  If $w$ is a weak $k$-partition and $\gamma$ is a weak $k$-composition, we shall say that 
\emph{$\gamma$ is of type $w$} if $\gamma$ is a permutation of $w$.
\end{definition}

\begin{definition}
A \emph{partition} $\lambda$ is any nonincreasing finite or infinite sequence 
$(\lambda_1, \lambda_2, \lambda_3, \dots)$ of nonnegative integers.  However, in contrast to Definition~\ref{defweakcomp}, only positive integers are considered parts, thus 
for a partition $\lambda$, $\ell(\lambda) = \#\{ i : \lambda_i > 0 \}.$  
Analogous to the frequency sequence of a weak composition, 
the \emph{frequency sequence} $f(\lambda)$ of $\lambda$ is
\[ f(\lambda) = (f_1(\lambda), f_2(\lambda), f_3(\lambda), \dots )  .\]
\end{definition}

\begin{remark}
In fact, no distinction will be drawn between, e.g., $\lambda=(5,2,1,1)$ and
$\lambda=(5,2,1,1,0,0,0,0,0,0,0,0,\dots)$; both will be considered the same partition of length $4$ and weight $9$.   
Also, $f\Big( (5,2,1,1) \Big) = (2,1,0,0,1,0,0,\dots)$.

It will be convenient to consider $\lambda_j = 0$ for any $j>\ell(\lambda)$, even when $\lambda$ is not explicitly constructed with a tail of zeros.
\end{remark}

\begin{definition}
 The set of all partitions of weight
 $n$ is denoted by $\mathscr{P}(n)$ and the cardinality of $\mathscr{P}(n)$ by $p(n)$.  The notation  $\lambda \vdash n$ means ``$\lambda$ is a partition of weight $n$", i.e., $\lambda\in\mathscr{P}(n)$.
\end{definition}

 For example, \[\mathscr{P}(4) = 
 \{ (4), (3,1), (2,2), (2,1,1), (1,1,1,1) \}, \] so $p(4)=5$.
\begin{definition}
  For a partition $\lambda$, the partition $\lambda - \one$ is the partition obtained from $\lambda$ by decreasing each of its parts by $1$:
  \[ \lambda-\one := ( \lambda_1 - 1, \lambda_2 -1, \dots, \lambda_{\ell(\lambda)}-1 ). \]
\end{definition}
Notice that $l(\lambda - \one) = l(\lambda) - f_1(\lambda).$
 
   It is often convenient to denote a partition (resp. weak $k$-partition)
by the \emph{superscript frequency notation} 
$\langle 1^{f_1} 2^{f_2} 3^{f_3} \cdots  \rangle$
(resp. $\langle 0^{f_0} 1^{f_1} 2^{f_2} 3^{f_3} \cdots  \rangle$) where it is 
permissible to omit $f_j$ if $f_j=1$ and to omit $j^{f_j}$ if $f_j=0$.
Thus, \[ 
(5,5,5,5,3,2,2,1,1,1) =
\langle 1^3 2^2 3\ 5^4 \rangle \] are two ways of expressing that particular (weak $10$-)partition of $30$.

  A variant on this notation for weak compositions (in order to emphasize runs of adjacent equal parts) will also be useful.  For example, let us allow ourselves to write the weak $9$-composition 
$(3,3,2,0,0,3,3,3,3)$ of $20$ as $[3^2 2^1 0^2 3^4]$.

The following quantities will arise often enough to warrant these definitions:
\begin{definition}
   
   Following Schneider~\cite{RS18}, for a partition $\lambda = (\lambda_1, \lambda_2, \dots, \lambda_k)$, we define its
   \emph{norm} to be the product of its parts,
   \[ \pp{\lambda}  := \lambda_1 \lambda_2 \cdots \lambda_k. \]  Further, the
 \emph{factorial} of a partition $\lambda$ is $\lambda! = (\lambda_1!, \lambda_2!, \dots, \lambda_k!)$, so 
 that $\pfp{\lambda} = \lambda_1! \lambda_2! \cdots \lambda_k!.$
 Analogously for a weak $k$-composition $\gamma = (\gamma_1, \dots, \gamma_k)$, \
     \[ \pfp{\gamma} := \gamma_1! \gamma_2! \dots \gamma_k! \]

 Observe that $f$ effectively maps a partition $\lambda = (\lambda_1, \lambda_2, \dots)$ 
 to a weak $\lambda_1$-composition of weight $l(\lambda)$ if we ignore the infinite tail
 of zeros in $f(\lambda)$.
 For example, $f( (5,5,5,5,3,2,2,1,1,1)) = (3,2,1,0,4)$, a weak $5$-composition of weight
 $10$.
    Likewise, $f$ 
maps a
 weak $k$-composition $\gamma$ to a weak $L$-composition of weight $k$, where $L$ is the
 largest part of $\gamma$.
       Thus the  \emph{frequency factorial product} of a weak $k$-composition $\gamma$
may be consistently       
       notated as
     \[ \mfp{\gamma} = f_0(\gamma)! f_1(\gamma)!  f_2(\gamma)! \cdots ,\]
     and that of a partition $\lambda$ as
     \[ \mfp{\lambda} = f_1(\lambda)! f_2(\lambda)! f_3(\lambda)! \cdots  . \] 

\end{definition}

\begin{observation} The number of weak $k$-compositions of type 
$w= (w_1,\dots, w_k)$, where $w$ is a weak $k$-partition, is $k!/\mfp{w}$.
\end{observation}

\begin{definition}
  A \emph{multipartition} is a $t$-tuple of partitions for some $t$. 
\end{definition}

For example, $( (1,1,1), (4,1), (3,2,2) )$ is a multipartition simply because $(1,1,1)$,
$(4,1)$, and $(3,2,2)$  are all partitions.

\begin{definition}
  The \emph{multipartition dissection} $\MD{\lambda}$ of the partition
  $\lambda$ is the following Cartesian product: 
 \[ \MD{\lambda} := \mathscr{P}(\lambda_1) \times \mathscr{P}(\lambda_2) \times
  \cdots \times \mathscr{P}(\lambda_{l(\lambda)}). \]
\end{definition}

We will require the result given by N. J. Fine~\cite[p. 38, Eq. (22.2)]{F},
\begin{equation} \label{fine}
   \sum_{\lambda\vdash n} \frac{1}{ 1^{f_1} 2^{f_2} 3^{f_3} \cdots 
   f_1 ! f_2! f_3! \cdots} = 1,
\end{equation}
 which may be 
expressed in the present notation as
\begin{equation}
  \sum_{\lambda\vdash n} \frac{1}{\pp{\lambda} \mfp{\lambda} } = 1,
\end{equation} in the iterated form
\begin{equation} \label{GenFine}
\sum_{i=1}^{p(\lambda_1)p(\lambda_2)\cdots p(\lambda_{l(\lambda)})}
\prod_{j=1}^{l(\lambda)} \frac{1}{\mfp{\mu^{ij}} \pp{\mu^{ij}} } = 1.
\end{equation}
The superscript notation on $\mu$ is to be understood as follows:
$\lambda = (\lambda_1, \dots, \lambda_l)$ is a partition; then $\mu^{ij}$ is the
$i$th partition of $\lambda_j$ where the $p(\lambda_1) p(\lambda_2) \cdots 
p(\lambda_l)$ multipartitions of $\MD{\lambda}$ have been placed in
\emph{some} order; any order is fine.
See also~\eqref{32term} below
for an explicit illustration.

Notice that~\eqref{fine} states that the sum of the co\"efficients that appear in 
the MacMahon partial fractions decomposition
\[ \frac{1}{(1-x)(1-x^2) \cdots (1-x^k)} = \sum_{\lambda\vdash k} 
\frac{1}{ 1^{f_1} 2^{f_2} 3^{f_3} \cdots 
   f_1 ! f_2! f_3! \cdots} g(\lambda; x, x, \dots, x),
 \] where $g$ is defined below in Eq.~\eqref{gdef}, must be $1$.

\subsubsection{Combinatorial Generating Functions}   
   As part of the combinatorial construction to be undertaken, we will need to associate with each partition $\lambda$ a certain rational generating function in indeterminates $x_1, x_2, \dots, x_{|\lambda|}$; namely let
   \begin{equation} \label{gdef}
   g(\lambda; \x):=  g(\lambda;  x_1, x_2, \dots, x_{|\lambda|}) :=
      \prod_{j=1}^{\ell(\lambda)}\frac{1}{1- \prod_{k=1}^{\lambda_j} x_{s(\lambda; j,k)}}, \end{equation}
 with \[ s(\lambda; j,k)  = k+\sum_{r=1}^{j-1} \lambda_r .\]

Of necessity, the notation used in defining~\eqref{gdef} for a general partition $\lambda$ makes a simple idea rather opaque.  To understand immediately how to construct $g(\lambda; \x)$ for any partition $\lambda$, simply consider, 
   for example, for the five partitions of $4$, we have the following associated ``$g$-functions":
   \begin{align*}
     g\Big((4); x_1, x_2, x_3, x_4\Big) &= \frac{1}{1-x_1 x_2 x_3 x_4} \\
     g\Big((3,1); x_1, x_2, x_3, x_4\Big) &= \frac{1}{(1-x_1 x_2 x_3) (1-x_4) } \\
     g\Big((2,2); x_1, x_2, x_3, x_4\Big) &= \frac{1}{(1-x_1 x_2)(1- x_3  x_4) } \\
     g\Big((2,1,1); x_1, x_2, x_3, x_4\Big) &= \frac{1}{(1-x_1 x_2)(1- x_3)(1 -x_4) } \\
       g\Big((1,1,1,1); x_1, x_2, x_3, x_4\Big) &= \frac{1}{(1-x_1 )(1-x_2)(1- x_3)(1 -x_4) } .
   \end{align*}
   
  We denote the symmetric group of degree $n$ by $\Sym_n$. 
 The application of a permutation $\sigma\in\Sym_{|\lambda|}$ to 
 $g(\lambda; \x)$, will be written as 
 $\sigma g(\lambda;\x)$, with the intended meaning
   \[ \sigma g(\lambda; \x) = g(\lambda; \sigma\x) = g( \lambda; x_{\sigma(1)}, x_{\sigma(2)}, \dots,
x_{\sigma(|\lambda|)}).\]

  Let $\Orb_{\lambda}$ 
  (resp. $H_\lambda$) 
 denote the orbit 
 (resp. stabilizer) 
 of $g(\lambda; \x)$ under the action of $\Sym_{|\lambda|}$.

  Thus,
  \begin{equation} \label{stabsize}
  \left| H_\lambda \right| = \pfp{\lambda} \mfp{\lambda},
  \end{equation}
 or, by the orbit--stabilizer theorem,
  
  \begin{equation}  \label{orbitsize}
   \left|  \Orb_{\lambda} \right | = \frac{ |\lambda|! }{  \pfp{\lambda} \mfp{\lambda}   } .
   \end{equation}
   
\begin{remark} \label{lambdaminusonermk}
  Notice that \[ \pfp{(\lambda-\one)}  |\Orb_{\lambda}|  
  = \frac{|\lambda|!}{ \pp{\lambda} \mfp{\lambda}  }, \] which is $|\lambda|!$ times the
  co\"efficient of the term indexed by $\lambda$ in the MacMahon decomposition.
\end{remark}

\subsection{Statement of main result}
  The goal is to understand~\eqref{Meq} combinatorially.   This will be accomplished by proving 
  the following natural multivariate generalization of~\eqref{Meq}:
  
  \begin{theorem} \label{MainResult}
\begin{multline} \label{GenMDecomp}
\sum_{\sigma\in\Sym_k} \sigma
\frac{1}{(1-x_1)(1-x_1 x_2)(1-x_1 x_2 x_3)\cdots(1-x_1 x_2 \cdots x_k)}\\
= \sum_{\lambda \vdash k} \pfp{(\lambda-\one)}
\sum_{\phi(\x)\in\Orb_\lambda } \phi(\x) 
\end{multline}
where 
$\Orb_\lambda$ is the orbit of $g(\lambda;\x)$ under the action of $\Sym_{|\lambda|}$,
and $g(\lambda; \x)$ is defined in~\eqref{gdef}.

\end{theorem}

\section{Partial Fractions Decompositions}
It is well known that for fixed positive integer $k$, the generating function for $p_k(n)$ is 
\begin{equation} F_k(x):=\sum_{n\geq 0} p_k (n) x^n = \prod_{j=1}^k \frac{1}{1-x^j}. \label{pmngf} \end{equation}
Since the right-hand side of~\eqref{pmngf} is a rational function, it can be decomposed into ordinary partial fractions, as considered, e.g, by Cayley~\cite{C} and Rademacher~\cite[p. 302]{R73}, or into 
$q$-partial fractions, as studied by Munagi~\cite{AM07,AM08}.

In examining the ordinary partial fraction decompositions of, say, $F_4(x)$,
\begin{align}
F_4(x) &= \frac{-17/72}{x-1} + \frac{59/288}{(x-1)^2} + \frac{1/8}{(x-1)^3} + \frac{1/24}{(x-1)^4} + \frac{1/8}{x+1} + \frac{1/32}{(x+1)^2}\notag\\ & \qquad\qquad+ \frac{(x+1)/9}{x^2+x+1} + \frac{1/8}{x^2+1}\\
   &= \frac{-17/72}{x-1} + \frac{59/288}{(x-1)^2} + \frac{1/8}{(x-1)^3} + \frac{1/24}{(x-1)^4} + \frac{1/8}{x+1} + \frac{1/32}{(x+1)^2}\notag \\
  &\qquad\qquad+ \frac{(2+\omega^2)/27 }{x-\omega} + \frac{(2+\omega)/27}{x-\omega^2} + \frac{-i/16}{x-i} + \frac{i/16}{x+i},
\end{align}
where $\omega:=\exp(2\pi i /3),$
we notice immediately the apparent arbitrariness of the co\"efficients that arise in the expansion.

  For Munagi's $q$-partial fractions, the co\"efficients are nicer, but still not transparent:
\begin{multline}
F_4(x) =  \frac{25/144}{(1-x)^2} + \frac{1/8}{(1-x)^3} + \frac{1/24}{(1-x)^4} + \frac{1/16}{1-x^2} + \frac{1/8}{(1-x^2)^2}\notag\\  \qquad\qquad+ \frac{(x+2)/9}{1-x^3} + \frac{1/4}{1-x^4}.
\end{multline}

MacMahon's partial fraction decomposition of $F_k(x)$, 
\begin{equation} \label{mpf2}
  \prod_{j=1}^k \frac{1}{1-x^j} = 
         \sum_{ \lambda \vdash m } \frac{g(\lambda; x,x,x,\dots,x  )}
         {\mfp{\lambda} \pp{\lambda} } ,
 \end{equation}
thus has the distinct advantage that the co\"efficients are known \emph{a priori}, and furthermore these co\"efficients are ``combinatorial numbers" in the sense that they
are products of integer exponential and factorial expressions.

In order to begin to understand~\eqref{mpf2} combinatorially, we shall
multiply both sides of~\eqref{mpf2} by $k !$ 
and observe that 
 $g(\lambda; x,x,x,\dots,x)$
 is the generating function for the sequence that counts
  a certain class of restricted weak $k$-compositions defined below.

Equation~\eqref{mpf2} together with~\eqref{orbitsize}, after some investigation, suggested the generalization of MacMahon's partial fraction decomposition presented 
above as Theorem~\ref{MainResult}.

\section{Proof of Theorem~\ref{MainResult}}
Starting with the left member of~\eqref{GenMDecomp}, we have
\begin{align*}
&\phantom{=}
\sum_{\sigma\in\Sym_k} \sigma \left( \frac{1}{(1-x_1)(1-x_1 x_2) (1-x_1 x_2 x_3) \cdots 
(1-x_1 x_2 x_3 \cdots x_k)}\right)\\
&= \sum_{\sigma\in\Sym_k} \sigma \left(  \sum_{a_1, a_2, \dots, a_k\geq 0} 
x_1^{a_1} (x_1 x_2)^{a_2} (x_1 x_2 x_3)^{a_3} \cdots (x_1 x_2 \cdots x_k)^{a_k} \right)\\
&= \sum_{\sigma\in\Sym_k} \sigma \left(  \sum_{a_1, a_2, \dots, a_k \geq 0} 
 x_1^{a_1+a_2+\cdots + a_k} x_2^ {a_2+a_3+\cdots +a_k} \cdots 
  x_{k-1}^{a_{k-1} + a_k} x_k^{a_k} \right)\\
&= \sum_{\sigma\in\Sym_k} \sigma \left(  \sum_{w_1\geq w_2 \geq \dots \geq w_k \geq 0} 
 x_1^{w_1} x_2^ {w_2} \cdots x_k^{w_k} \right) \\
 &= \sum_{\sigma\in\Sym_k} \sigma \left( \sum_{w\in\mathscr{W}_k}  x_1^{w_1} x_2^ {w_2} \cdots x_k^{w_k} \right) \\
  &=\sum_{\gamma\in\mathscr{C}_k} \mfp{\gamma}   x_1^{\gamma_1} x_2^{\gamma_2} \cdots x_k^{\gamma_k}.
\end{align*}

Thus, we see that the left member of \eqref{GenMDecomp} generates every 
weak $k$-composition (where the $j$th part appears as the exponent of $x_j$) exactly 
$\mfp{\gamma}$ times.

Now let us consider the right member of~\eqref{GenMDecomp}, 
\begin{equation} \label{GenMDecompRHS}
 \sum_{\lambda \vdash m} 
 \pfp{  (\lambda - \one) }
 \sum_{\phi(\x)\in\Orb_\lambda}
\phi(\x) , \end{equation} where $\Orb_{\lambda}$ is the orbit of $g(\lambda; \x)$ under the (transitive) action of $\Sym_{|\lambda|}$.

Pick an arbitrary weak $k$-composition $\gamma$.  We need to show that the term 
$x_1^{\gamma_1} x_2^{\gamma_2} \cdots x_k^{\gamma_k}$ appears in the expansion of~\eqref{GenMDecompRHS} with co\"efficient 
$\mfp{\gamma}$.
Associated with $\gamma$ is the frequency sequence 
$f(\gamma)= (f_0(\gamma), f_1(\gamma), f_2(\gamma), \dots )$.  Permute the nonzero terms of $f(\gamma)$ into nondecreasing order to form a partition $\lambda$ of weight $k$, and we may write $\lambda = \lambda(\gamma)$, since the partition $\lambda$ is uniquely determined by $\gamma$.
Thus it must be the case that there exists $\sigma\in\Sym_k$ such that the 
weak $k$-composition $\sigma(\gamma)$ is of type 
$[ c_1^{\lambda_1} c_2^{\lambda_2} \cdots c_\ell^{\lambda_\ell} ]$ for some distinct nonnegative integers $c_1, c_2, \dots, c_\ell$.

 For a given $\lambda\vdash k$ of length $\ell$, we have, by expanding~\eqref{gdef} as a series,
\begin{multline*}
 g(\lambda; \x) = \sum_{c_1, c_2, \dots, c_\ell\geq 0} (x_1 x_2 \cdots x_{\lambda_1})^{c_1} (x_{\lambda_1+1} x_{\lambda_1+2} \cdots
 x_{\lambda_1+\lambda_2})^{c_2} \cdots \\ \times
 (x_{\lambda_1+\lambda_2 + \cdots \lambda_{l-1} +1} x_{\lambda_1 + \lambda_2 + \cdots + \lambda_{\ell-1} + 2} \cdots x_{\lambda_1 + \lambda_2 + \cdots + \lambda_{\ell}})^{c_\ell},
 \end{multline*}
so $g(\lambda; \x)$ is the generating function for weak $k$-compositions of type
\[ [ c_1^{\lambda_1} c_2^{\lambda_2} \cdots c_\ell^{\lambda_\ell} ]. \]

 Now the orbit $\Orb_\lambda$ of $g(\lambda; \x)$ under the action of $\Sym_k$ 
 contains the terms that generate all permutations of weak $k$-compositions of 
 type $[ c_1^{\lambda_1} c_2^{\lambda_2} \cdots c_\ell^{\lambda_\ell} ]$.
 
  The terms $x_1^{\gamma_1} x_2^{\gamma_2} \cdots x_k^{\gamma_k}$ are generated by 
  those terms of~\eqref{GenMDecompRHS}  
 in the orbit of $g(\mu;\x)$ 
for multipartitions $\mu \in \MD{\lambda}$ 
where $\lambda= \lambda(\gamma)$.
 
   For each weak $k$-composition $\gamma$, and corresponding partition 
   $\lambda=
  \lambda(\gamma) $ $= (\lambda_1, \lambda_2, \dots,\lambda_{l(\lambda)} )$, we
  generate all the associated multipartitions in $\MD{\lambda}$.  
  Let $\mu^{ij}_k$ denote the $k$th part in the partition $\mu^{ij}$, where $\mu^{ij}$
is the $i$th
 partition of $\lambda_j$, the $j$th part of $\lambda$.  
 
   Note that $i$ runs from
  $1$ through $p(\lambda_1)p(\lambda_2) \cdots p(\lambda_{l(\lambda)})$ where some
  ordering has been imposed on the multipartitions (any ordering will do).  Of course,
  $j$ runs from $1$ to $l(\lambda)$, and $k$ runs from $1$ to $l(\mu^{ij})$. 
  
  For example, if we wish to calculate the number of times the weak
  $5$-composition $\gamma=(7,7,4,7,4)$ is generated by the right-hand side 
  of~\eqref{GenMDecomp}, i.e., the number of times the 
  expression
 $x_1^7 x_2^7 x_3^4 x_4^7 x_5^4$ appears, 
 we see, by symmetry, 
 that this must be the same as the number of times $(7,7,7,4,4)$ appears.

The weak $5$-partition $(7,7,7,4,4)$ is clearly obtained from $(7,7,4,7,4)$ by allowing the permutation
$\sigma = (3,4)$ to act on it.
  Then $\lambda(\gamma) = \lambda(\sigma\gamma) = (3,2)$ because $\gamma$ 
  (and $\sigma\gamma$) both contain
  $3$ of one part and $2$ of another part.  
  
    We notice that a certain number of copies of $(7,7,7,4,4)$ are generated by each
 of the terms 
 \begin{multline} \label{32xs}
  \frac{2}{(1-x_1 x_2 x_3)(1-x_4 x_5)}, 
  \frac{2}{(1-x_1 x_2 x_3)(1-x_4)(1- x_5)},\\
  \frac{1}{(1-x_1 x_2)(1- x_3)(1-x_4 x_5)}, \
  \frac{1}{(1-x_1 x_2) (1-x_3)(1-x_4) (1- x_5)},\\
  \frac{1}{(1-x_1)(1- x_2 )(1-x_3)(1-x_4 x_5)},\\
  \frac{1}{(1-x_1)(1- x_2)(1- x_3)(1-x_4)(1- x_5)},
 \end{multline} and by no other terms.
To aid our analysis we consider the multipartition dissection of 
the partition $\lambda=(3,2)$:
  \begin{multline} \label{32term} \MD{ (3,2)} = \{\mbox{all partitions of $3$} \} \times
  \{ \mbox{all partitions of $2$} \} \\ =
\Big\{ \Big( (3), (2) \Big),  \Big( (3), (1,1) \Big), 
 \Big( (2,1), (2) \Big),  \Big( (2,1), (1,1) \Big), \\
  \Big( (1,1,1), (2) \Big),  \Big( (1,1,1), (1,1) \Big) \Big\}, \end{multline}
  because each of these six multipartitions indexes a term that generates some number
of copies of $\sigma\gamma=(7,7,7,4,4).$
  (In this example, we have $\mu^1 = ((3),(2))$, $\mu^2 = ((3),(1,1))$, 
  $\mu^3 = ((2,1),(2))$, $\mu^4 = ((2,1),(1,1))$, $\mu^5 = ((1,1,1),(2))$, and 
  $\mu^6 = ((1,1,1),(1,1))$;  
  $\mu^{11} = \mu^{21} = (3)$, $\mu^{12} = \mu^{32} = \mu^{52} = (2)$,
  $\mu^{22} = \mu^{42} = \mu^{62} = (1,1)$, and
   $\mu^{51} = \mu^{61} = (1,1,1)$.)
   
    We use elementary combinatorial reasoning to count how many copies 
 of $(7,7,7,4,4)$ are generated by each of the six terms.  That number is a consequence
 of the commutivity of ordinary multiplication.  
For example, consider the third listed term in~\eqref{32xs}

 To generate $(7,7,7,4,4)$,
we may do so by any of the following permutations of this third term:
\begin{multline*}  \frac{1}{(1-x_1 x_2)(1- x_3)(1-x_4 x_5)},  
\frac{1}{(1-x_1 x_3)(1- x_2)(1-x_4 x_5)}, \\
\frac{1}{(1-x_2 x_3)(1- x_1)(1-x_4 x_5)}, \end{multline*}
which are indexed by the multipartition $\mu^3 = ((2,1),(2))$.

 This clearly lists all elements in the Cartesian product of the two orbits: one is the
 orbit of $\frac{1}{(1-x_1x_2)(1-x_3)}$ under the action of $\Sym_3 = \Sym_{\{ 1,2,3 \} }$,
(the permutations of $\{1,2,3\}$), and the other is the orbit of $\frac{1}{1-x_4 x_5}$ under the action of $\Sym_{\{ 4,5\}}$
(the permutations of $\{4,5\}$).
Since each term generates one copy of $(7,7,7,4,4)$, the total contribution of these terms is given by
\begin{multline}
  (\mu^{31}-\one)! (\mu^{32} - \one)! \left| \Orb_{\mu^{31}} \right|  \left| \Orb_{\mu^{32}} \right|  \\ = 
  \frac{\lambda_1!}{  f_1(\mu^{31})! f_2(\mu^{31})! \cdots \mu^{31}_1 \mu^{31}_2    }  \cdot 
     \frac{ \lambda_2!}{ f_1(\mu^{32})! f_2(\mu^{32})! \cdots \mu^{32}_1   },
\end{multline} where we have applied Remark~\ref{lambdaminusonermk}.

  Of course, to generate \emph{all} copies of $(7,7,7,4,4)$, we must sum over all of the terms indexed by the six
  members of $\MD{(3,2)}$, employing the analogous counting formula in each case.

   In the general case, the preceding combinatorial argument yields
  \[  \prod_{i,j} (\mu^{ij}-\one)! \left| \Orb_{\mu^{ij}} \right| 
  =  \sum_{i=1}^{p(\lambda_1)p(\lambda_2) \cdots p(\lambda_{l(\lambda)})} 
 \prod_{j=1}^{l(\lambda)} \frac{\lambda_j!}{\mfp{ \mu^{ij}} \pp{\mu^{ij}}}. \]

 Thus all that remains in order to prove Theorem~\ref{MainResult} is to establish:
  \begin{equation} \label{almost}
   \mfp{\gamma} = \sum_i \prod_j \frac{\lambda_j!}{\mfp{ \mu^{ij}} \pp{\mu^{ij}}} .
  \end{equation}

  Since $\mfp{\gamma} = \pfp{\lambda} = \prod_j \lambda_j!$, we immediately see that~\eqref{almost} is
equivalent to the assertion
\[  1 = \sum_i \prod_j \frac{1}{\mfp{ \mu^{ij}} \pp{\mu^{ij}}} , \]
  which is exactly~\eqref{GenFine}, and thus Theorem~\ref{MainResult} is 
  established. \qed
  
  \section{Example: the case $k=4$}
  Before concluding, let us examine the $k=4$ case in some detail.  Our main result,
Theorem~\ref{MainResult}, in the case $k=4$ asserts
 \begin{multline} \label{MacMahon4}
 \sum_{\sigma\in\Sym_4} \sigma\left( \frac{1}{(1-x_1)(1-x_1 x_2)(1-x_1 x_2 x_3)(1-x_1 x_2 x_3 x_4)}
 \right) \\ = \frac{6}{1-x_1 x_2 x_3 x_4} 
 + 2\left( \frac{1}{(1-x_1 x_2 x_3)(1-x_4)}   + \frac{1}{(1-x_1 x_2 x_4)(1-x_3)} \right. \\
 +  \left. \frac{1}{(1-x_1 x_3 x_4)(1-x_2)} + \frac{1}{(1-x_1)(1-x_2 x_3 x_4)} \right)\\
 + \left( \frac{1}{(1-x_1 x_2)(1-x_3 x_4)} +  \frac{1}{(1-x_1 x_3)(1-x_2 x_4)} \right. \\
 + \left. \frac{1}{(1-x_1 x_4)(1-x_2 x_3)}  \right) \\
 + \left( \frac{1}{(1-x_1 x_2)(1-x_3)(1-x_4)} + \frac{1}{(1-x_1 x_3)(1-x_2)(1-x_4)}\right. \\
 + \left. \frac{1}{(1-x_1 x_4)(1-x_2)(1-x_3)} + \frac{1}{(1-x_1)(1-x_2 x_3)(1-x_4)}\right. \\
  + \left. \frac{1}{(1-x_1)(1- x_2 x_4)(1-x_3)} + \frac{1}{(1-x_1)(1-x_2) (1-x_3 x_4)}\right) \\
  + \frac{1}{ (1-x_1)(1-x_2)(1-x_3)(1-x_4) }.
 \end{multline}
In the left member of~\eqref{MacMahon4}, we have 
\[ \frac{1}{(1-x_1)(1-x_1 x_2)(1-x_1 x_2 x_3)(1-x_1 x_2 x_3 x_4)}, \] which generates every weak $4$-partition $w=(w_1, w_2, w_3, w_4)$ exactly once.  
The cardinality of the orbit of the action of $\Sym_4$ on $w$ is 
 \[ \frac{4!}{ f_0(w)! f_1(w)! f_2(w)! \cdots },\] i.e., there are  $4!/ (f_0(w)! f_1(w)! f_2(w)! \cdots )$ distinct weak $4$-compositions of type $w$.  Or equivalently, a given weak $4$-compositon $\gamma$ which equals $\sigma w$ for some 
permutation $\sigma\in\Sym_4$, is generated 
$\mfp{\gamma} = f_0(\gamma)! f_1(\gamma)! f_2(\gamma)! \cdots  $  times.

  The generation of weak $4$-compositions on the right side of~\eqref{MacMahon4} is more subtle.  Notice that the terms of the right side are grouped according to the partitions of $4$ (which index the sum on the right side) 
  in the order $(4)$, $(3,1)$, $(2,2)$, $(2,1,1)$, $(1,1,1,1)$. 
 For a given weak $4$-composition $\gamma$, the multiplicities of the parts determine which of the terms of the right side contribute to its generation. 
 
  A detailed summary is provided in Table~\ref{table}.  In order to make sure the table is clear, let us look at one
 line of it in detail.  Observe the case with $\lambda=(2,2)$ and
 form of $\gamma$ as $abab$.  The $abab$ means we are considering weak 
 $4$-compositions
 where the first and third parts are the same, and the second and fourth parts are the
 same, but the first and second parts are different.  The corresponding terms from the 
 right member of~\eqref{MacMahon4} are equivalent to
 \begin{equation} \label{MT4}
 (23)g(22,\x) + (23)g(211,\x) + (14)g(211,\x) + g(1111,\x).
 \end{equation}
 The first term of~\eqref{MT4} is $(23)g(22,\x)$, i.e. apply the transposition $(23)$ to
 \[ g(22,\x) = \frac{1}{(1-x_1 x_2) (1-x_3 x_4)}; \] this yields
 \begin{equation} \label{MT41} (23) g(22,\x) =  \frac{1}{(1-x_1 x_3) (1-x_2 x_4)}.  
 \end{equation}
 Expand each factor of the right side of~\eqref{MT41} as a geometric series to find that
 weak $4$-compositions $(w_1,w_2,w_3,w_4)$ are generated (in the exponents of the 
 $x_i$'s) in which $w_1 = w_3$ and $w_2 = w_4$, i.e., compositions of the type 
 $(a,b,a,b)$.  Is this the only way that compositions of type $(a,b,a,b)$ may be generated?
 No.  Consider the second term of~\eqref{MT4}, $(23) g(211,\x)$, which is
 \[  \frac{1}{(1-x_1 x_3)(1-x_2)(1-x_4)} . \]
 This term generates weak $4$-compositions $(w_1,w_2,w_3,w_4)$, in which $w_1=w_3$.
 \emph{Some} of the weak compositions generated by this term will happen to have
 $w_2 = w_4$, and thus these will be of the form $(a,b,a,b)$ as well, i.e., this term 
generates compositions of the general form $(a,b,a,c)$; on those occasions that it
happens to be the case that $b=c$, we have a weak composition of the form considered
by this particular line of the table.  And so on, with the third and forth terms of
~\eqref{MT4}.  The other lines of the table may be interpreted similarly.

 \begin{table}[h!] 
\begin{center}\begin{tabular}{|c | c | l | l |}
\hline\hline
$\lambda$ & Form of $\gamma$ & Generating terms of RHS of~\eqref{MacMahon4} \\
\hline\hline
 & $aaaa$ & $6 g( 4; \x) + 2\Big( ()+(34) +(24)+(14) \Big)g(31;\x)$ \\
 $(4)$     &             &  $+ \Big( () + (23) + (24)  \Big) g(22; \x) $\\
      &              &  $+\Big( () + (13) + (23) + (24)+(14)+(13)(24)  \Big) g(211; \x) $\\
      & & $+ g(1111; \x)$ \\              
\hline
        & $aaab$ & $2 g(31; \x) +\Big( ()+(13)+(23) \Big) g(211; \x) + g(1111;\x) $\\
$(31)$ & $aaba$ & $2 (34) g(31;\x) + \Big( ()+ (24) + (14) \Big) g(211;\x) + g(1111;\x)$\\
       & $abaa$ & $2 (24) g(31; \x) + \Big( (23) + (24) + (13)(24)   \Big) g(211;\x) + g(1111;\x)$ \\
       & $abbb$ & $2 (14) g(31;\x) + \Big( (13) + (14) + (13)(24) \Big) g(211;\x) + g(1111;\x)$ \\
  \hline
     & $aabb$ & $g(22;\x) + \Big(  () + (13)(24)  \Big)g(211;\x) +g(1111;\x)$ \\
 $(22)$ & $abab$ & $(23) g(22;\x) + \Big(  (23) + (14)  \Big)g(211;\x) + g(1111;\x)$ \\     
   &  $abba$ & $(24) g(22;\x) + \Big( (13) + (24) \Big)g(211;\x) + g(1111;\x)$ \\
   \hline 
   & $aabc$ & $g(211;\x) + g(1111;\x)$ \\  
   & $abac$ & $(23) g(211;\x) + g(1111;\x)$ \\
 $(211)$ & $abca$ & $(24) g(211;\x) + g(1111;\x)$ \\  
   & $abbc$ & $(13) g(211;\x) + g(1111;\x)$ \\
     & $abcb$ & $(14) g(211;\x) + g(1111;\x)$ \\
       & $abcc$ & $(13)(24) g(211;\x) + g(1111;\x)$ \\ \hline
 $(1111)$ & $abcd$ & $g(1111;\x)$ \\ \hline\hline      
\end{tabular}
\caption{The letters $a, b, c, d$ represent distinct nonnegative integers.  Permutations are presented in cycle notation. } \label{table}
\end{center}
\end{table}

\section*{Acknowledgments}
The author thanks George Andrews for pointing out~\cite[p. 209, Ex. 1]{A76}, which lead
to the research culminating in this paper.
The author thanks Matthew Katz for his interest 
and useful suggestions.  The author \emph{particularly} thanks Robert Schneider
for discussions and encouragement of this project over a long period of time, and
for carefully reading and offering concrete suggestions to improve earlier 
versions of the manuscript.  Finally, the author is \emph{extremely grateful} to the editor and 
anonymous referees for carefully
reading the manuscript, catching errors, offering numerous helpful suggestions, and for
their kind patience as the author prepared revisions.

\end{document}